\documentclass[12pt]{article}

\usepackage{theorem,amssymb}

\author{Jean-Paul Allouche \\
CNRS, LRI, UMR 8623 \\
Universit\'e Paris-Sud, B\^atiment 490 \\
F-91405 Orsay Cedex (France) \\
{\tt allouche@lri.fr} \\
\and
Christiane Frougny \\ 
LIAFA, CNRS UMR 7089 \\
Case 7014 \\ 
75205 Paris Cedex 13 (France) \\
and Universit\'e Paris 8 \\
{\tt Christiane.Frougny@liafa.jussieu.fr}   
}

\title{Univoque numbers and an avatar of Thue-Morse}

\date{ }

\def \proof{\bigbreak\noindent{\it Proof.\ \ }}

\def \endpf{{\ \ $\Box$ \medbreak}}

\newtheorem{theorem}{Theorem}

\newtheorem{corollary}{Corollary}
\newtheorem{proposition}{Proposition}
\theorembodyfont{\rm}

\newtheorem{remark}{Remark}

\newtheorem{definition}{Definition}

\begin{document}

\baselineskip=17pt

\maketitle

\begin{abstract}
Univoque numbers are real numbers $\lambda > 1$ such that the number $1$ 
admits a unique expansion in base $\lambda$, i.e., a unique expansion 
$1 = \sum_{j \geq 0} a_j \lambda^{-(j+1)}$, with 
$a_j \in \{0, 1, \ldots, \lceil \lambda \rceil -1\}$ for every $j \geq 0$. 
A variation of this definition was studied in 2002 by Komornik and Loreti, 
together with sequences called {\em admissible sequences}. We show how a 1983 
study of the first author gives both a result of Komornik and Loreti on the 
smallest admissible sequence on the set $\{0, 1, \ldots, b\}$, and a result 
of de Vries and Komornik (2007) on the smallest univoque number belonging to 
the interval $(b, b+1)$, where $b$ is any positive integer. We also prove that 
this last number is transcendental. An avatar of the Thue-Morse sequence, 
namely the fixed point beginning in $3$ of the morphism $3 \to 31$, $2 \to 30$, 
$1 \to 03$, $0 \to 02$, occurs in a ``universal'' manner. 

\medskip

\noindent
{\it 2000 Math.~Subject Classification.} 11A63, 11B83, 11B85, 68R15, 11J81.

\noindent
{\it Keywords.} Beta-expansion, univoque numbers, iteration of continuous 
functions, Thue-Morse sequence, uniform morphism, automatic sequence, 
transcendence.

\end{abstract}

\section{Introduction}

Komornik and Loreti determined in \cite{KL1} the smallest {\em univoque\,}
real number in the interval $(1, 2)$, i.e., the smallest number 
$\lambda \in (1, 2)$ such that $1$ has a unique expansion 
$1 = \sum_{j \geq 0} a_j / \lambda^{j+1}$   
with $a_j \in \{0, 1\}$ for every $j \geq 0$. The word ``univoque'' in this
context seems to have been introduced (with a slightly different meaning) by 
Dar\'oczy and K\'atai in \cite{DK1, DK2}, while characterizing unique expansions
of the real number $1$ was done by Erd\H{o}s, Jo\'o, and Komornik in 
\cite{EJK}. The first author and Cosnard showed in \cite{AC3} how the result 
of \cite{KL1} parallels (and can be deduced from) 
their study of a certain set of binary sequences arising in the iteration of 
unimodal continous functions of the unit interval that was done in 
\cite{Cosnard, AC1, All}. The relevant sets of binary sequences occurring in 
references \cite{AC1, All}, resp. in reference \cite{KL1}, can be defined by
$$
\begin{array}{lll}
\Gamma &:=& \{A \in \{0, 1\}^{\mathbb N}, \ \forall k \geq 0, \ 
\overline{A} \leq \sigma^k A \leq A\} \\
& & \\
\Gamma_{strict} &:=& \{A \in \{0, 1\}^{\mathbb N}, \ \forall k \geq 1, 
\ \overline{A} < \sigma^k A < A\} \\
\end{array}
$$
where $\sigma$ is the {\em shift} on sequences and the {\em bar} operation replaces 
$0$'s by $1$'s and $1$'s by $0$'s, i.e., if $A = (A_n)_{n \geq 0}$, 
then $\sigma A := (a_{n+1})_{n \geq 0}$, and $\overline{A} := (1-a_n)_{n \geq 0}$; 
furthermore $\leq$ denotes the lexicographical order on sequences induced by 
$0 < 1$, the notation $A < B$ meaning as usual that $A \leq B$ and $A \neq B$. 
The smallest univoque number in $(1, 2)$ and the smallest nonperiodic sequence of 
the set $\Gamma$ both involve the Thue-Morse sequence (see for example \cite{AS1}
for more on this sequence).

\bigskip

It is tempting to generalize these sets to alphabets with more than $2$ letters. 

\begin{definition}

For $b$ a positive integer, we will say that the real number $\lambda > 1$ is
{\em $\{0, 1, \ldots, b \}$-univoque} if the number $1$ has a unique expansion 
as $1 = \sum_{j \geq 0} a_j \lambda^{-(j+1)}$, where $a_j$ belongs to 
$\{0, \ldots, b\}$ for all $j \geq 0$. Furthermore, if $\lambda > 1$ is 
$\{0, 1, \ldots, \lceil \lambda \rceil - 1 \}$-univoque, we will simply say 
that $\lambda$ is {\em univoque}.
\end{definition}

\begin{remark}
If $\lambda > 1$ is $\{0, 1, \ldots, b \}$-univoque for some positive 
integer $b$, then $\lambda \leq b+1$. Also note that any integer
$q \geq 2$ is univoque, since there is exactly one expansion of $1$ as
$1 = \sum_{j \geq 0} a_j q^{-(j+1)}$, with $a_j \in \{0, 1, \ldots, q-1\}$, 
namely $1 = \sum_{j \geq 0} (q-1) q^{-(j+1)}$.
\end{remark}

Komornik and Loreti studied in \cite{KL2} the reals $\lambda$ belonging to 
the interval $(1, b+1]$ that are $\{0, 1, \ldots, b \}$-univoque. For their study,
they introduced {\em admissible sequences\,} on the alphabet $\{0, 1, \ldots, b\}$. 
Denote, as above, by $\sigma$ the {\em shift\,} on sequences, and by {\em bar\,}
the operation that replaces every $t \in \{0, 1, \ldots, b\}$ by $b-t$, i.e., if 
$A = (a_n)_{n \geq 0}$, then $\overline{A} := (b-a_n)_{n \geq 0}$.
Also denote by $\leq$ the lexicographical order on sequences induced by the 
natural order on $\{0, 1, \ldots, b\}$. Then, a sequence $A = (a_n)_{n \geq 0}$ 
on $\{0, 1, \ldots, b\}$ is {\em admissible\,} if
$$
\begin{array}{llll}
\forall k \geq 0 \ \mbox{\rm such that } a_k < b, \ 
                             &\sigma^{k+1} A  &<& A, \\
\forall k \geq 0 \ \mbox{\rm such that } a_k > 0, \ 
                             &\sigma^{k+1} A &>& \overline{A}. \\
\end{array}
$$
(Note that our notation is not exactly the notation of \cite{KL2}, since our
sequences are indexed by ${\mathbb N}$ and not ${\mathbb N} \setminus \{0\}$.)
These sequences have the following property:
the map that associates with a real $\lambda \in (1, b+1]$ 
the sequence of coefficients $(a_j)_{j \geq 0} \in \{0, 1, \ldots, b\}$ 
of the greedy (i.e., the lexicographically largest) expansion of $1$, 
$1 = \sum_{j \geq 0} a_j \lambda^{-(j+1)}$, is a bijection from the set of 
$\{0, 1, \ldots, b \}$-univoque $\lambda$'s to the set of admissible sequences 
on $\{0, 1, \ldots, b\}$ (see \cite{KL2}).

\bigskip

Now there are two possible generalizations of the result of \cite{KL1} about 
the smallest univoque number in $(1, 2)$, i.e., the smallest admissible binary
sequence. One is to look at the smallest (if any) admissible sequence on 
the alphabet $\{0, 1, \ldots, b\}$, as did Komornik and Loreti in \cite{KL2},
the other is to look at the smallest (if any) univoque number in $(b, b+1)$,
as did de Vries and Komornik in \cite{VriKom}. 

\bigskip

It happens that the first author already studied a generalization of the set 
$\Gamma$ to the case of more than $2$ letters (see \cite[Part~3]{All}).
Interestingly enough this study was not related to the iteration of
continuous functions as was the study of $\Gamma$, but only introduced as a 
tempting formal arithmetico-combinatorial generalization of the study of the set 
of binary sequences $\Gamma$ to a similar set of sequences with more than two 
values.

\bigskip

The purpose of the present paper is threefold:

1) {\em to show how the 1983 study \cite[Part~3, p. 63--90]{All}
gives both the result of Komornik and Loreti in \cite{KL2} on the 
smallest admissible sequence on $\{0, 1, \ldots, b\}$, and the result
of de Vries and Komornik in \cite{VriKom} on the smallest number
univoque number $\lambda$ belonging to $(b, b+1)$ where $b$ is 
any positive integer;} 

2) {\em to bring to light a {\em universal} morphism that governs the
smallest elements in 1) above, and to show that the infinite sequence 
generated by this morphism is an avatar of the Thue-Morse sequence;}

3) {\em to prove that the smallest univoque number belonging to $(b, b+1)$ 
(where $b$ is any positive integer) is transcendental.}

\bigskip

The paper consists of five sections. In Section~2 below we recall some results 
of \cite[Part~3, p. 63--90]{All} on the generalization of the set $\Gamma$ to 
a $(b+1)$-letter alphabet. Then we give some properties of the lexicographically 
least nonperiodic sequence of this set, completing results of 
\cite[Part~3, p. 63--90]{All}.
In Section~3 we give two corollaries of the properties of this least 
sequence: one gives the result in \cite{KL2}, the other gives the result
in \cite{VriKom}. The transcendence results are proven in the last section.

\bigskip

\section{The generalized $\Gamma$ and $\Gamma_{strict}$ sets}

\begin{definition}\label{fund-def}
Let $b$ be a positive integer, and $\cal A$ be a finite ordered set
with $b+1$ elements. Let $\alpha_0 < \alpha_1 < \ldots < \alpha_b$ be 
the elements of $\cal A$. The {\em bar} operation is defined on $\cal A$ 
by $\overline{\alpha_j} = \alpha_{b-j}$. 
We extend this operation to ${\cal A}^{\mathbb N}$
by $\overline{(a_n)_{n \geq 0}} := (\overline{a_n})_{n \geq 0}$.
Let $\sigma$ be the {\em shift} on ${\cal A}^{\mathbb N}$, defined by
$\sigma((a_n)_{n \geq 0}) := (a_{n+1})_{n \geq 0}$.

We define the sets $\Gamma({\cal A})$ and 
$\Gamma_{strict}({\cal A})$ by:
$$
\begin{array}{llll}
\Gamma({\cal A}) &:=& \! \!
\{A = (a_n)_{n \geq 0} \in {\cal A}^{\mathbb N}, \ 
a_0 = \max{\cal A}, \ \forall k \geq 0, \ \overline{A} \leq \sigma^k A \leq A\}, \\
\\
\Gamma_{strict}({\cal A}) \! \!
&:=& \! \!
\{A = (a_n)_{n \geq 0} \in {\cal A}^{\mathbb N}, \ 
a_0 = \max{\cal A}, \ \forall k \geq 1, \ \overline{A} < \sigma^k A < A\}. \\
\end{array}
$$
\end{definition}

\begin{remark}
The set $\Gamma({\cal A})$ was introduced by the first author in 
\cite[Part~3, p. 63--90]{All}. Note that there is a misprint in the 
definition given on p.~66 in \cite{All}: $a_{\beta - i}$ should be changed
into $a_{\beta - 1 - i}$ as confirmed by the rest of the text.
\end{remark}

\begin{remark}\label{nonperiodic}
A sequence belongs to $\Gamma_{strict}({\cal A})$ if and only if 
it belongs to $\Gamma({\cal A})$ and is nonperiodic. Namely, 
$\sigma^k A = A$ if and only if $A$ is $k$-periodic; 
if $\sigma^k A = \overline{A}$, then $\sigma^{2k} A = A$, and the
sequence $A$ is $2k$-periodic. 
\end{remark}

\begin{remark}\label{precision}
If the set ${\cal A}$ is given by ${\cal A} := \{i, i+1, \ldots, i+z\}$
where $i$ and $z$ are integers, equipped with the natural order, then for 
any $x \in {\cal A}$, we have $\overline{x} = 2i + z - x$. Namely, following 
Definition~\ref{fund-def} above, we write $\alpha_0 := i, \alpha_1 := i+1, 
\ldots, \alpha_z := i+z$. Hence, for any $j \in [0,z]$, we have 
$\overline{\alpha_j} = \alpha_{z-j}$, which can be rewritten 
$\overline{i+j} = i+z-j$, i.e., for any $x$ in ${\cal A}$, we have
$\overline{x} = i+z-(x-i) = 2i+z-x$.
\end{remark}

A first result is that the sets $\Gamma_{strict}({\cal A})$ are closely 
linked to the set of admissible sequences whose definition was recalled 
in the introduction.

\begin{proposition}\label{gamma-admissible}
Let $A = (a_n)_{n \geq 0}$ be a sequence in 
$\{0, 1, \ldots, b\}^{\mathbb N}$, such that $a_0 = t \in [0,b]$.
Suppose that the sequence $A$ is not equal to $b \ b \ b \ \ldots$ 
Then the sequence $A$ is admissible if and only if $2t > b$ and $A$
belongs to the set $\Gamma_{strict}(\{b-t, b-t+1, \ldots, t\})$. 
(The order on $\{b-t, b-t+1, \ldots, t\}$ is induced by the order on 
${\mathbb N}$. From Remark~\ref{precision} the bar operation is given by 
$\overline{j} = b-j$.) 
\end{proposition}

\proof Let $A = (a_n)_{n \geq 0}$ be a sequence belonging to
$\{0, 1, \ldots, b\}^{\mathbb N}$ such that $a_0 = t \in [0, b-1]$,
and such that $A \neq b \ b \ b \ldots$ 

\medskip

$*$ First suppose that $2t > b$ and $A$ belongs to 
$\Gamma_{strict}(\{b-t, b-t+1, \ldots, t\})$. Then, for all 
$k \geq 1$, $\overline{A} < \sigma^k A < A$, which clearly 
implies that $A$ is admissible.

\medskip

$*$ Now suppose that $A$ is admissible. We thus have
$$
\begin{array}{llll}
\forall k \geq 1 \ \mbox{\rm such that } a_{k-1} < b, \
                             &\sigma^k A  &<& A, \\
\forall k \geq 1 \ \mbox{\rm such that } a_{k-1} > 0, \
                             &\sigma^k A &>& \overline{A}. \\
\end{array}
$$
We first prove that, if the sequence $A$ is not a constant sequence, then 
$$
\forall k \geq 1, \ \overline{A} < \sigma^k A < A.
$$ 
We only prove the inequalities $\sigma^k A < A$. The remaining inequalities
are proved in a similar way. If $a_{k-1} < b$, the inequality $\sigma^k A < A$ 
holds. If $a_{k-1} = b$, there are two cases: 
\begin{itemize}

\item either $a_0 = a_1 = \ldots = a_{k-1} = b$, then, if $a_k < b$ we clearly
have $\sigma^k A < A$; if $a_k = b$, then the sequence $\sigma^k A$ begins with
some block of $b$'s followed by a letter $< b$, thus it begins with a block of
$b$'s shorter than the initial block of $b$'s of the sequence $A$ itself, hence 
$\sigma^k A < A$;

\item or there exists an index $\ell$ with $1 < \ell < k$, such that
$a_{\ell-1} < b$, and $a_{\ell} = a_{\ell + 1} = \ldots = a_{k-1} = b$. 
As $A$ is admissible, we have $\sigma^{\ell} A < A$. It thus suffices to
prove that $\sigma^k A \leq \sigma^{\ell} A$. This is clearly the case
if $a_k < b$. On the other hand, if $a_k = b$, the sequence $\sigma^k A$
begins with a block of $b$'s which is shorter than the initial block of
$b$'s of the sequence $\sigma^{\ell} A$, hence $\sigma^k A \leq \sigma^{\ell} A$. 
\end{itemize}
Now, since $a_0 = t$ and $\sigma^k A < A$ for all $k \geq 1$, we have
$a_k \leq t$ for all $k \geq 0$. Similarly, since $\sigma^k A > \overline{A}$
for all $k \geq 1$, we have $a_k \geq b - t$ for all $k \geq 1$. Finally
$A > \overline{A}$ implies that $t = a_0 \geq b - t$. We thus have that
$2t \geq b$ and $A$ belongs to $\Gamma(\{b - t, b - t + 1, \ldots, t\})$. 
Now, if $b=2t$, then $\{b-t, b-t+1, ..., t \} = \{t\}$ and $\bar{t} = t$.
This implies that $A = t \ t \ t \ \ldots$, which is not an admissible
sequence. \endpf 

\begin{remark}
For $b=1$, this (easy) result is noted without proof in \cite{EJK} 
and proved in \cite{AC3}.
\end{remark}

We need another definition from \cite{All}. 

\begin{definition}
Let $b$ be a positive integer, and $\cal A$ be a finite ordered set  with 
$b+1$ elements. Let $\alpha_0 < \alpha_1 < \ldots < \alpha_b$ be the elements 
of $\cal A$. We suppose that ${\cal A}$ is equipped with a bar operation as 
in Definition~\ref{fund-def}. Let $A = (a_n)_{n \geq 0}$ be a periodic
sequence of {\em smallest} period $T$, and such $a_{T-1} < \max{\cal A}$. 
Let $a_{T-1} = \alpha_j$ (thus $j < b$). Then the sequence $\Phi(A)$ is 
defined as the $2T$-periodic sequence beginning with 
$a_0 \ a_1 \ \ldots \ a_{T-2} \ \alpha_{j+1} \ \overline{a_0} \ \overline{a_1} \ 
\ldots \ \overline{a_{T-2}} \ \alpha_{b - j - 1}$, i.e.,
$$
\Phi((a_0 \ a_1 \ \ldots \ a_{T-2} \ \alpha_j)^{\infty} :=
(a_0 \ a_1 \ \ldots \ a_{T-2} \ \alpha_{j+1} \ \overline{a_0} \ \overline{a_1} \     
\ldots \ \overline{a_{T-2}} \ \alpha_{b - j - 1})^{\infty}.
$$
\end{definition}

We first prove the following easy claim.

\begin{proposition}\label{smallest}
The smallest element of $\Gamma(\{b - t, b - t + 1, \ldots, t\})$ 
(where $2t > b$) is the $2$-periodic sequence $(t \ (b-t))^{\infty} = 
(t \ (b-t) \ t \ (b-t) \ t \ \ldots)$.
\end{proposition}

\proof Since any sequence $A = (a_n)_{n \geq 0}$ belonging to
$\Gamma(\{b - t, b - t + 1, \ldots, t\})$ begins in $t$, and satisfies
$\sigma A \geq \overline{A}$, then it must satisfy $a_0 = t$ and
$a_1 \geq b-t$. Now if a sequence $A$ belonging to
$\Gamma(\{b - t, b - t + 1, \ldots, t\})$ is such that $a_0 = t$ and
$a_1 = (b-t)$, then it must be equal to the $2$-periodic sequence 
$(t \ (b-t))^{\infty}$ (\cite[Lemma~2, b, p.~73]{All}). Since this periodic
sequence trivially belongs to $\Gamma(\{b - t, b - t + 1, \ldots, t\})$,
it is its smallest element. \endpf

Denoting as usual by $\Phi^s$ the $s$-th iterate of $\Phi$, we state the 
following theorem which is a particular case of the theorem on pages 
72--73 of \cite{All} about the smallest elements in certain subintervals 
of $\Gamma(\{0, 1, \ldots, b\})$, and of the definition of $q$-mirror 
sequences given in \cite[Section~II, 1, p.~67]{All} (here $q=2$).

\begin{theorem}[\cite{All}]\label{old-thm}
Define $P := (t \ (b-t))^{\infty} = (t \ (b-t) \ t \ (b-t) \ t \ \ldots)$.
The smallest nonperiodic sequence in the set
$\Gamma(\{b - t, b - t + 1, \ldots, t\})$ (i.e., the smallest element of 
$\Gamma_{strict}(\{b - t, b - t + 1, \ldots, t\})$) is the sequence $M$ defined by
$$
M := \lim_{s \to \infty} \Phi^s(P),
$$ 
that actually takes the (not 
necessarily distinct) values $b-t$, $b-t+1$, $t-1$, $t$. Furthermore, 
this sequence $M = (m_n)_{n \geq 0} = 
t \ \ \ b-t+1 \ \ \ b-t \ \ \ t \ \ \ b-t \ \ \ t-1 \ \ldots$ 
can be recursively defined by 
$$
\begin{array}{lll}
\forall k \geq 0, \ m_{2^{2k}-1} = t, \\
\forall k \geq 0, \ m_{2^{2k+1}-1} = b+1-t, \\ 
\forall k \geq 0, \ \forall j \in [0, 2^{k+1}-2], \
m_{2^{k+1}+j} = \overline{m_j}. \\
\end{array}
$$
\end{theorem}

It was proven in \cite{All} that the sequence 
$\lim_{s \to \infty} \Phi^s((t \ (b-t))^{\infty}$ is $2$-automatic
(for more about automatic sequences, see \cite{AS2}). The second author
noted that this sequence is actually a fixed point of a uniform morphism of 
length $2$ as soon as the cardinality of the set $\{b-t, b-t+1, \ldots, b\}$ 
is at least equal to $4$, i.e., $2t \geq b+3$. (Recall that we always have
$t \geq b-t$, i.e., $2t \geq b$.) More precisely we have Theorem~\ref{new-thm}
below, where the Thue-Morse sequence pops up, as in \cite{All} and in \cite{KL2}, 
but also as in \cite{AC1} and \cite{KL1}. Before stating this theorem we give
a definition.

\begin{definition}\label{universal}
The ``universal'' morphism $\Theta$ is defined on $\{e_0, e_1, e_2, e_3\}$ by
$$
\Theta(e_3) := e_3 e_1, \ \Theta(e_2) := e_3 e_0, \
\Theta(e_1) := e_0 e_3, \ \Theta(e_0) := e_0 e_2.
$$
Note that this morphism has an infinite fixed point beginning in $e_3$
$$
\Theta^{\infty}(e_3) = \lim_{k \to \infty} \Theta^k(e_3)
= e_3 \ e_1 \ e_0 \ e_3 \ e_0 \ e_2 \ e_3 \ e_1 \ e_0 \ e_2 \ldots.
$$
\end{definition}

\begin{theorem}\label{new-thm}
Let $(\varepsilon_n)_{n \geq 0}$ be the Thue-Morse sequence, defined by 
$\varepsilon_0 = 0$ and for all $k \geq 0$, $\varepsilon_{2k} = \varepsilon_k$
and $\varepsilon_{2k+1} = 1 - \varepsilon_k$.
Then the smallest nonperiodic sequence $M = (m_n)_{n \geq 0}$ belonging to 
$\Gamma(\{b - t, b - t + 1, \ldots, t\})$ satisfies
$$
\forall n \geq 0, \ m_n = \varepsilon_{n+1} - (2t - b -1) \varepsilon_n + t - 1.
$$
Using the morphism $\Theta$ introduced in Definition~\ref{universal}
above we thus have

\begin{itemize}
\item if $2t \geq b+3$, then the sequence $M$ is the fixed point beginning
      in $t$ of the morphism deduced from $\Theta$ by renaming 
      $e_0, e_1, e_2, e_3$ respectively $b-t, b-t+1, t-1, t$ (note that the 
      condition $2t \geq b+3$ implies that these four numbers are distinct); 
\item if $2t = b+2$ (thus $b-t+1 = t-1$), then the sequence $M$ is the 
      pointwise image of the fixed point beginning in $e_3$ of the morphism 
      $\Theta$ by the map $g$ defined by $g(e_3) := t$, $g(e_2) = g(e_1) := t-1$, 
      $g(e_0) := b-t$; 
\item if $2t = b+1$ (thus $b-t = t-1$ and $b-t+1 = t$), then the sequence $M$ 
      is the pointwise image of the fixed point beginning in $e_3$ of the 
      morphism $\Theta$ by the map $h$ defined by $h(e_3) = h(e_1) := t$, 
      $h(e_2) = h(e_0) := t-1$.
\end{itemize}
\end{theorem}

\proof Let us first prove that the sequence $M = (m_n)_{n \geq 0}$ is equal to
the sequence $(u_n)_{n \geq 0}$, where
$u_n := \varepsilon_{n+1} - (2t-b-1)\varepsilon_n + t-1$. It suffices to
prove that the sequence $(u_n)_{n \geq 0}$ satisfies the recursive relations
defining $(m_n)_{n \geq 0}$ that are given in Theorem~\ref{old-thm}.
Recall that the sequence $(\varepsilon_n)_{n \geq 0}$ has the property
that $\varepsilon_n$ is equal to the parity of the sum of the binary digits
of the integer $n$ (see \cite{AS1} for example). Hence, for all $k \geq 0$,
$\varepsilon_{2^{2k}-1} = 0$, $\varepsilon_{2^{2k+1}-1} = 1$, and
$\varepsilon_{2^{2k}} = \varepsilon_{2^{2k+1}} = 1$. This implies that
for all $k \geq 0$, $u_{2^{2k}-1} = t$ and $u_{2^{2k+1}-1} = b+1-t$. 
Furthermore, for all $k \geq 0$, and for all $j \in [0, 2^{k+1}-2]$, we have
$\varepsilon_{2^{k+1}+j} = 1 - \varepsilon_j$ and $\varepsilon_{2^{k+1}+j+1}
= 1 - \varepsilon_{j+1}$. Hence $u_{2^{k+1}+j} = b - u_j = \overline{u_j}$.

\medskip

To show how the ``universal'' morphism $\Theta$ enters the picture, 
we study the sequence $(v_n)_{n \geq 0}$ with values in $\{0, 1\}^2$
defined by: for all $n \geq 0$, $v_n := (\varepsilon_n, \varepsilon_{n+1})$.
Since we have, for all $n \geq 0$, $v_{2n} = (\varepsilon_n, 1-\varepsilon_n)$
and $v_{2n+1} = (1-\varepsilon_n, \varepsilon_{n+1})$, we clearly have
$$
\begin{array}{llll}
\mbox{\rm if} &v_n = (0, 0), \, \mbox{\rm then} 
&v_{2n} = (0, 1) \ \mbox{\rm and} &v_{2n+1} = (1, 0), \\
\mbox{\rm if} &v_n = (0, 1), \, \mbox{\rm then}
&v_{2n} = (0, 1) \ \mbox{\rm and} &v_{2n+1} = (1, 1), \\
\mbox{\rm if} &v_n = (1, 0), \, \mbox{\rm then}
&v_{2n} = (1, 0) \ \mbox{\rm and} &v_{2n+1} = (0, 0), \\
\mbox{\rm if} &v_n = (1, 1), \, \mbox{\rm then}
&v_{2n} = (1, 0) \ \mbox{\rm and} &v_{2n+1} = (0, 1). \\
\end{array}
$$
This exactly means that the sequence $(v_n)_{n \geq 0}$ is the
fixed point beginning in $(0, 1)$ of the $2$-morphism
$$
\begin{array}{lll}
(0, 0) &\to& (0, 1) (1, 0) \\
(0, 1) &\to& (0, 1) (1, 1) \\
(1, 0) &\to& (1, 0) (0, 0) \\
(1, 1) &\to& (1, 0) (0, 1) \\
\end{array}
$$
We may define $e_0 := (1, 0)$, $e_1 := (1, 1)$, $e_2 := (0, 0)$, $e_3 := (0, 1)$.
Then the above morphism can be written 
$$
e_3 \to e_3 e_1, \ e_2 \to e_3 e_0, \ e_1 \to e_0 e_3, \ e_0 \to e_0 e_2
$$ 
which is the morphism $\Theta$. The above construction shows that 
the sequence $(v_n)_{n \geq 0}$ is a fixed point of $\Theta$.

\medskip

Now, define the map $\omega$ on $\{0, 1\}^2$ by 
$$
\omega((x,y)) := y - (2t-b-1)x + t - 1.
$$
We have $\omega(v_n) = m_n$ for all $n \geq 0$. Thus
\begin{itemize}
\item if $2t \geq b+3$, the sequence $(m_n)_{n \geq 0}$ takes exactly
      four distinct values, namely $b-t, b-t+1, t-1, t$. This implies
      that $(m_n)_{n \geq 0}$ is the fixed point beginning in $t$ of
      the morphism obtained from $\Theta$ by renaming the letters, i.e.,
      $e_3 \to t$, $e_2 \to (t-1)$, $e_1 \to (b-t+1)$, $e_0 \to (b-t)$. The
      morphism can thus be writen $t \to t \ (b-t+1)$, $(t-1) \to t \ (b-t)$, 
      $(b-t+1) \to (b-t) \ t$, $(b-t) \to (b-t) \ (t-1)$;

\item if $2t = b+2$ (resp. $2t = b+1$) the sequence $(m_n)_{n \geq 0}$
      takes exactly three (resp. two) values, namely $b-t, t-1, t$
      (resp. $t-1, t$). It is still the pointwise image by $\Theta$ of 
      the sequence $(v_n)_{n \geq 0}$. Renaming $\Theta$ as $g$ (resp. $h$)
      as in the statement of Theorem~\ref{new-thm} only takes into account 
      that the integers $b-t, b-t+1, t-1, t$ are not distinct. \endpf
\end{itemize}

\begin{remark}
The reason for the choice of indexes for $e_3, e_2, e_1, e_0$ is that
the order of indexes is the same as the natural order on the integers
$t, t-1, b-t+1, b-t$ to which they correspond when $2t \geq b+3$.
In particular if $b=t=3$, the morphism reads:
$3 \to 31$, $2 \to 30$, $1 \to 03$, $0 \to 02$.
Interestingly enough, though not surprisingly, this morphism also occurs
(up to renaming once more the letters) in the study of infinite square-free
sequences on a $3$-letter alphabet. Namely, in the paper \cite{Berstel}, 
Berstel proves that the square-free Istrail sequence \cite{Istrail}, originally
defined (with no mention of the Thue-Morse sequence) as the fixed point of the 
(non-uniform) morphism $0 \to 12$, $1 \to 102$, $2 \to 0$, is actually the 
pointwise image of the fixed point beginning in $1$ of a $2$-morphism $\Theta'$ 
on the $4$-letter alphabet $\{0, 1, 2, 3\}$ by the map $0 \to 0$, $1 \to 1$,
$2 \to 2$, $3 \to 0$. The morphism $\Theta'$ is given by
$$
\Theta'(0) = 12, \ \Theta'(1) = 13, \ \Theta'(2) = 20, \ \Theta'(3) = 21.
$$
The reader will note immediately that $\Theta'$ is another avatar of $\Theta$
obtained by renaming letters as follows: $0 \to 2$, $1 \to 3$, $2 \to 0$,
$3 \to 1$. This, in particular, shows that {\em the sequence $(m_n)_{n \geq 0}$,
in the case where $2t = b+2$, is the fixed point of the non-uniform morphism
$t \to t \ (t-1) \ (b-t)$, $(t-1) \to t \ (b-t)$, $(b-t) \to (t-1)$, i.e., 
an avatar of Istrail's square-free sequence. Furthermore it results from 
\cite{Berstel} that this sequence on three letters cannot be the fixed point 
of a uniform morphism.} A last remark is that the square-free Braunholtz sequence
on three letters given in \cite{Braunholtz} (see also \cite[p. 18-07]{Berstel})
is exactly our sequence $(m_n)_{n \geq 0}$ when $t=b=2$, i.e., the sequence
$2 \ 1 \ 0 \ 2 \ 0 \ 1 \ 2 \ 1 \ 0 \ 1 \ 2 \ 0 \ \ldots$
\end{remark}

\section{Small admissible sequences and small univoque numbers 
with given integer part}

\subsection{Small admissible sequences with values in the set
$\{0, 1, \ldots, b\}$}

In \cite{KL2} the authors are interested in the smallest admissible
sequence with values in the set $\{0, 1, \ldots, b\}$, where $b$ is an
integer $\geq 1$. They prove in particular the following result, which
is an immediate corollary of our Theorem~\ref{new-thm}.

\begin{corollary}[Theorems~4.3 and 5.1 of \cite{KL2}]
Let $b$ be an integer $\geq 1$. The smallest admissible sequence with values 
in $\{0, 1, \ldots, b\}$ is the sequence $(z+\varepsilon_{n+1})_{n \geq 0}$ 
if $b = 2z+1$, and $(z+\varepsilon_{n+1}-\varepsilon_n)_{n \geq 0}$ if $b=2z$.
\end{corollary}

\proof Let $A = (a_n)_{n \geq 0}$ be the smallest (non-constant) admissible 
sequence with values in $\{0, 1, \ldots, b\}$. Since $A > \overline{A}$, we
must have $a_0 \geq \overline{a_0} = b - a_0$. 

Thus, if $b=2z+1$ we have $a_0 \geq z+1$. We also have, for all $i \geq 0$, 
$\overline{a_0} \leq a_i \leq a_0$. Now the smallest element of the set 
$\Gamma(\{b-z-1, b-z, \ldots, z-1, z+1\})$ is the smallest admissible sequence 
on $\{0, 1, \ldots, b\}$ that begins in $z+1$. Hence this is the smallest 
admissible sequence with values in $\{0, 1, \ldots, b\}$. Theorem~\ref{new-thm}
gives that this sequence is $(m_n)_{n \geq 0}$ with, for all $n \geq 0$,
$m_n = \varepsilon_{n+1} +z$.

If $b=2z$, we have $a_0\geq z$. But if $a_0=z$, then $\overline{a_0}=z$, and
the conditions of admissibility implies that $a_n = z$ for all $n \geq 0$ and
$(a_n)_{n \geq 0}$ would be the constant sequence $(z \ z \ z \ldots)$. Hence
we must have $a_0 \geq z+1$. Now the smallest element of the set
$\Gamma(\{b-z-1, b-z, \ldots, z-1, z+1\})$ is the smallest admissible sequence   
on $\{0, 1, \ldots, b\}$ that begins in $z+1$. Hence this is the smallest
admissible sequence with values in $\{0, 1, \ldots, b\}$. Theorem~\ref{new-thm}
gives that this sequence is $(m_n)_{n \geq 0}$ with, for all $n \geq 0$,
$m_n = \varepsilon_{n+1} - \varepsilon_n + z$. \endpf

\subsection{Small univoque numbers with given integer part}

We are interested here in the univoque numbers $\lambda$ in an interval 
$(b, b+1]$ with $b$ a positive integer. This set was studied in \cite{KK}, 
where it was proven of Lebesgue measure $0$.
Since $1 = \sum_{j \geq 0} a_j \lambda^{-(j+1)}$, $\lambda \in (b, b+1]$ and 
$a_0 \leq b$, the fact that the expansion of $1$ is unique, hence equal to the 
greedy expansion, implies that $a_0 = b$. In other words, we study the admissible 
sequences with values in $\{0, 1, \ldots, b\}$ that begin in $b$, i.e., the set 
$\Gamma_{strict}(\{0, 1, \ldots, b\})$. We prove here, as a corollary of
Theorem~\ref{new-thm}, that, for any positive integer $b$, there exists a 
smallest univoque number belonging to $(b, b+1]$. 
This result was obtained in \cite{VriKom} (see the penultimate remark in that
paper); it generalizes the result obtained for $b=1$ in \cite{KL1}.

\begin{corollary}
For any positive integer $b$, there exists a smallest univoque number 
in the interval $(b, b+1]$. This number is the solution of the equation
$1 = \sum_{n \geq 0} d_n \lambda^{-n-1}$, where the sequence 
$(d_n)_{n \geq 0}$ is given by, for all $n \geq 0$, 
$d_n := \varepsilon_{n+1} - (b-1)\varepsilon_n + b - 1$.
\end{corollary}

\proof It suffices to apply Theorem~\ref{new-thm} with $t=b$. \endpf

\section{Transcendence results}

We prove here, mimicking the proof given in \cite{AC2}, that 
numbers such that the expansion of $1$ is given by
the sequence $(m_n)_{n \geq 0}$ are transcendental. This generalizes
the transcendence results of \cite{AC2} and \cite{KL2}.

\begin{theorem}
Let $b$ be an integer $\geq 1$ and $t \in [0, b]$ be an integer such that
$2t \geq b+1$. Define the sequence $(m_n)_{n \geq 0}$ as in Theorem~\ref{new-thm}
by, for all $n \geq 0$,  
$m_n := \varepsilon_{n+1} - (2t-b-1) \varepsilon_n + t - 1$, thus 
the sequence $(m_n)_{n \geq 0}$ begins with 
$t \ \ \ b-t+1 \ \ \ b-t \ \ \ t \ \ \ b-t \ \ \ t-1 \ \ldots$
Then the number $\lambda$ belonging to $(1, b+1)$ defined by
$1 = \sum_{n \geq 0} m_n \lambda^{-n-1}$ is transcendental. 
\end{theorem}

\proof Define the $\pm 1$ Thue-Morse sequence $(r_n)$ by 
$r_n := (-1)^{\varepsilon_n}$. We clearly have $r_n = 1 - 2 \varepsilon_n$
(recall that $\varepsilon_n$ is $0$ or $1$). It is also immediate that the
function $F$ defined for the complex numbers $X$ such that $|X| < 1$ by
$F(X) = \sum_{n \geq 0} r_n X^n$ satisfies 
$F(X) = \prod_{k \geq 0} (1 - X^{2^k})$ (see, e.g., \cite{AS1}). Since
$$
2m_n = 2\varepsilon_{n+1} - 2(2t-b-1) \varepsilon_n + 2t - 2 =
b - r_{n+1} + (2t-b-1)r_n 
$$
we have, for $|X| < 1$,
$$
2X\sum_{n \geq 0} m_n X^n = ((2t-b-1)X-1)F(X) + 1 + \frac{bX}{1-X}\cdot
$$
Taking $X = 1/\lambda$ where $1 = \sum_{n \geq 0} m_n \lambda^{-n-1}$, we
get the equation
$$
2 = ((2t-b-1)\lambda^{-1}-1)F(1/\lambda) + 1 + \frac{b}{\lambda - 1}\cdot
$$
Now, if $\lambda$ were algebraic, then this equation shows that $F(1/\lambda)$
would be an algebraic number. But, since $1/\lambda$ would be an algebraic
number in $(0, 1)$, the quantity $F(1/\lambda)$ would be transcendental from 
a result of Mahler \cite{Mahler}, giving a contradiction. \endpf

\begin{remark}
In particular the $\{0, 1, \ldots, b\}$-univoque number corresponding to 
the smallest admissible sequence with values in $\{0, 1, \ldots, b\}$ is 
transcendental, as proved in \cite{KL2} (Theorems~4.3 and 5.9). Also the 
smallest univoque number belonging to $(b, b+1)$ is transcendental.
\end{remark}

\section{Conclusion}

There are many papers dealing with univoque numbers. We will just mention 
here the study of univoque Pisot numbers. The authors together with K.~G.~Hare 
determined in \cite{AFH} the smallest univoque Pisot number, which happens 
to have algebraic degree $14$. Note that the number corresponding to the 
sequence of Proposition~\ref{smallest} is the larger real root of the 
polynomial $X^2 - tX - (b-t+1)$, hence a Pisot number (which is unitary 
if $t=b$). Also note that for any $b \geq 2$, the real number $\beta$ such
that the $\beta$-expansion of $1$ is $b 1^{\infty}$ is a univoque Pisot 
number belonging to the interval $(b, b+1)$. It would be interesting to
determine the smallest univoque Pisot number belonging to $(b, b+1)$: 
the case $b=1$ was addressed in \cite{AFH}, but the proof uses heavily the fine
structure of Pisot numbers in the interval $(1, 2)$ (see \cite{Am, Ta1, Ta2}). 
A similar study of Pisot numbers in $(b, b+1)$ would certainly help.

\bigskip

\noindent
{\bf Acknowledgments.} The authors thank M. de Vries and V. Komornik for
their remarks on a previous version of this paper.

\end{document}